# DISCUSSION OF "LEAST ANGLE REGRESSION" BY EFRON ET AL.

By Hemant Ishwaran

*Cleveland Clinic Foundation*

Being able to reliably, and automatically, select variables in linear regression models is a notoriously difficult problem. This research attacks this question head on, introducing not only a computationally efficient algorithm and method, LARS (and its derivatives), but at the same time introducing comprehensive theory explaining the intricate details of the procedure as well as theory to guide its practical implementation. This is a fascinating paper and I commend the authors for this important work.

Automatic variable selection, the main theme of this paper, has many goals. So before embarking upon a discussion of the paper it is important to first sit down and clearly identify what the objectives are. The authors make it clear in their introduction that, while often the goal in variable selection is to select a "good" linear model, where goodness is measured in terms of prediction accuracy performance, it is also important at the same time to choose models which lean toward the parsimonious side. So here the goals are pretty clear: we want good prediction error performance but also simpler models. These are certainly reasonable objectives and quite justifiable in many scientific settings. At the same, however, one should recognize the difficulty of the task, as the two goals, low prediction error and smaller models, can be diametrically opposed. By this I mean that certainly from an oracle point of view it is true that minimizing prediction error will identify the true model, and thus, by going after prediction error (in a perfect world), we will also get smaller models by default. However, in practice, what happens is that small gains in prediction error often translate into larger models and less dimension reduction. So as procedures get better at reducing prediction error, they can also get worse at picking out variables accurately.

Unfortunately, I have some misgivings that LARS might be falling into this trap. Mostly my concern is fueled by the fact that Mallows' $C_p$ is the criterion used for determining the optimal LARS model. The use of $C_p$ often leads to overfitting, and this coupled with the fact that LARS is a forward optimization procedure, which is often found to be greedy, raises some potential flags. This, by the way, does not necessarily mean that LARS per se is overfitting, but rather that I think $C_p$ may be an inappropriate model selection criterion for LARS. It is this point that will be the focus of my discussion. I will offer some evidence that $C_p$ can sometimes be used effectively if *model uncertainty* is accounted for, thus pointing to ways for its more appropriate use within LARS. Mostly I will make my arguments by way of high-dimensional simulations. My focus on high dimensions is motivated in part by the increasing interest in such problems, but also because it is in such problems that performance breakdowns become magnified and are more easily identified. Note that throughout my discussion I will talk only about LARS, but, given the connections outlined in the paper, the results should also naturally apply to the Lasso and Stagewise derivatives.

**1. Is $C_p$ the correct stopping rule for LARS?** The $C_p$ criterion was introduced by Mallows (1973) to be used with the OLS as an unbiased estimator for the model error. However, it is important to keep in mind that it was not intended to be used *when the model is selected by the data* as this can lead to selection bias and in some cases poor subset selection [Breiman (1992)]. Thus, choosing the model with lowest $C_p$ value is only a heuristic technique with sometimes bad performance. Indeed, ultimately, this leads to an inconsistent procedure for the OLS [Shao (1993)]. Therefore, while I think it is reasonable to assume that the $C_p$ formula (4.10) is correct [i.e., that it is reasonable to expect that $df(\widehat{\boldsymbol{\mu}}_k) \approx k$ under a wide variety of settings], there is really no reason to expect that minimizing the $C_p$ value will lead to an optimal procedure for LARS.

In fact, using $C_p$ in a Forward Stagewise procedure of any kind seems to me to be a risky thing to do given that $C_p$ often overfits and that Stagewise procedures are typically greedy. Figure 5 of the paper is introduced (partly) to dispel these types of concerns about LARS being greedy. The message there is that $\text{pe}(\widehat{\boldsymbol{\mu}})$, a performance measurement related to prediction error, declines slowly from its maximum value for LARS compared to the quick drop seen with standard forward stepwise regression. Thus, LARS acts differently than well-known greedy algorithms and so we should not be worried. However, I see the message quite differently. If the maximum proportion explained for LARS is roughly the same over a large range of steps, and hence models of different dimension, then this implies that there is not much to distinguish between higher-and lower-dimensional models. Combine this with the use of $C_p$ which could provide poor estimates for the prediction error due to selection bias and there is real concern for estimating models that are too large.

To study this issue, let me start by reanalyzing the diabetes data (which was the basis for generating Figure 5). In this analysis I will compare LARS to a Bayesian method developed in Ishwaran and Rao (2000), referred to as SVS (short for Stochastic Variable Selection). The SVS procedure is a hybrid of the spike-and-slab model approach pioneered by Mitchell and Beauchamp (1988) and later developed in George and McCulloch (1993). Details for SVS can be found in Ishwaran and Rao (2000, 2003). My reason for using SVS as a comparison procedure is that, like LARS, its coefficient estimates are derived via shrinkage. However, unlike LARS, these estimates are based on model averaging *in combination* with shrinkage. The use of model averaging is a way of accounting for model uncertainty, and my argument will be that models selected via $C_p$ based on SVS coefficients will be more stable than those found using LARS thanks to the extra benefit of model averaging.

Figures 1 and 2 present the $C_p$ values for the main effects model and the quadratic model from both procedures (the analysis for LARS was based on S-PLUS code kindly provided by Trevor Hastie). The $C_p$ values for SVS were computed by (a) finding the posterior mean values for coefficients, (b) ranking covariates by the size of their absolute posterior mean coefficient values (with the top rank going to the largest absolute mean) and (c) computing the $C_p$ value $C_p(\widetilde{\boldsymbol{\mu}}_k) = \|\mathbf{y} - \widetilde{\boldsymbol{\mu}}_k\|/\overline{\sigma}^2 - n + 2k$, where $\widetilde{\boldsymbol{\mu}}_k$ is the OLS estimate based on the $k$ top ranked covariates. All covariates were standardized. This technique of using $C_p$ with SVS was discussed in Ishwaran and Rao (2000).

We immediately see some differences in the figures. In Figure 1, the final model selected by SVS had $k = 6$ variables, while LARS had $k = 7$ variables. More interesting, though, are the discrepancies for the quadratic model seen in Figure 2. Here the optimal SVS model had $k = 8$ variables in contrast to the much higher $k = 15$ variables found by LARS. The top eight variables from SVS (some of these can be read off the top of the plot) are bmi, ltg, map, hdl, sex, age.sex, bmi.map and glu.2. The last three variables are interaction effects and a squared main effects term. The top eight variables from LARS are bmi, ltg, map, hdl, bmi.map, age.sex, glu.2 and bmi.2. Although there is a reasonable overlap in variables, there is still enough of a discrepancy to be concerned. The different model sizes are also cause for concern. Another worrisome aspect for LARS seen in Figure 2 is that its $C_p$ values remain bounded away from zero. This should be compared to the $C_p$ values for SVS, which attain a near-zero mininum value, as we would hope for.

**2. High-dimensional simulations.** Of course, since we do not know the true answer in the diabetes example, we cannot definitively assess if the LARS models are too large. Instead, it will be helpful to look at some simulations for a more systematic study. The simulations I used were designed following the recipe given in Breiman (1992). Data was simulated in all cases by using i.i.d. N(0,1) variables for $\varepsilon_i$. Covariates $x_i$, for $i = 1, \ldots, n$, were generated independently from a multivariate normal distribution with zero mean and with covariance satisfying $E(x_{i,j}x_{i,k}) = \rho^{|j-k|}$. I considered two settings for $\rho$: (i) $\rho = 0$ (uncorrelated); (ii) $\rho = 0.90$ (correlated). In all simulations, $n = 800$ and $m = 400$. Nonzero $\beta_j$ coefficients were in 15 clusters of 7 adjacent variables centered at every 25th variable. For example, for the variables clustered around the 25th variable, the coefficient values were given by $\beta_{25+j} = |h - j|^{1.25}$ for $|j| < h$, where $h = 4$. The other 14 clusters were defined similarly. All other coefficients were set to zero. This gave a total of 105 nonzero values and 295 zero values. Coefficient values were adjusted by multiplying by a common constant to make the theoretical $R^2$ value equal to 0.75 [see Breiman (1992) for a discussion of this point]. Please note that, while the various parameters chosen for the simulations might appear specific, I also experimented with other simulations (not reported) by considering different configurations for the dimension $m$, sample size $n$, correlation $\rho$ and the number of nonzero coefficients. What I found was consistent with the results presented here.

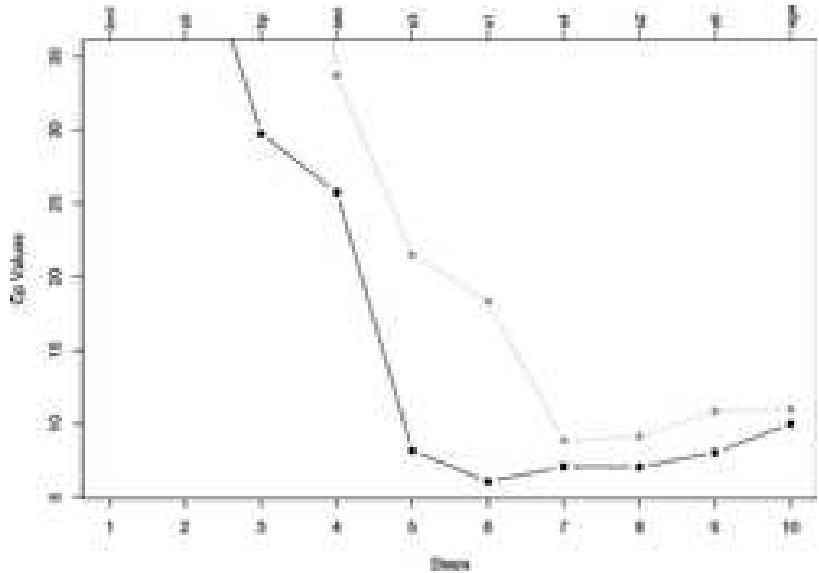

Fig. 1. $C_p$ values from main effects model for diabetes data: thick line is values from SVS; thin dashed line is from LARS. Covariates listed at the top of the graph are ordered by importance as measured by their absolute posterior mean.

For each $\rho$ correlation setting, simulations were repeated 100 times independently. Results are recorded in Table 1. There I have recorded what I call TotalMiss, FDR and FNR. TotalMiss is the total number of misclassified variables, that is, the total number of falsely identified nonzero $\beta_j$ coefficients and falsely identified zero coefficients; FDR and FNR are the false discovery and false nondiscovery rates defined as the false positive and false negative rates for those coefficients identified as nonzero and zero, respectively. The TotalMiss, FDR and FNR values reported are the averaged values from the 100 simulations. Also recorded in the table is $\widehat{m}$, the average number of variables selected by a procedure, as well as the performance value pe($\widehat{\boldsymbol{\mu}}$) [cf. (3.17)], again averaged over the 100 simulations.

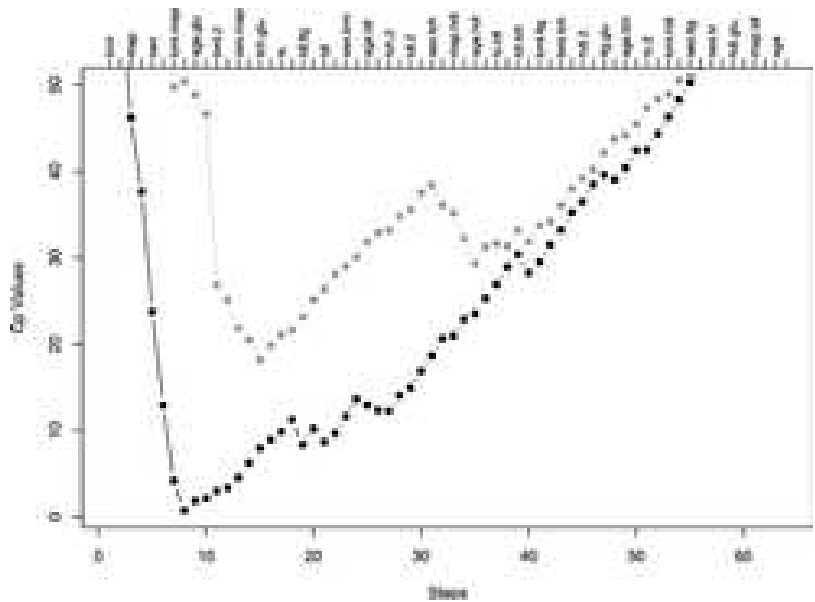

Fig. 2. *$C_p$ values from quadratic model: best model from SVS is $k = 8$ (thick line) compared with $k = 15$ from LARS (thin dashed line). Note how the minimum value for SVS is nearly zero.*

Table 1 records the results from various procedures. The entry "svsCp" refers to the $C_p$-based SVS method used earlier; "Step" is standard forward stepwise regression using the $C_p$ criterion; "svsBMA" is the Bayesian model averaged estimator from SVS. My only reason for including svsBMA is to gauge the prediction error performance of the other procedures. Its variable selection performance is not of interest. Pure Bayesian model averaging leads to improved prediction, but because it does no dimension reduction at all it cannot be considered as a serious candidate for selecting variables.

The overall conclusions from Table 1 are summarized as follows:

1. The total number of misclassified coefficients and FDR values is high in the uncorrelated case for LARS and high in the correlated case for stepwise regression. Their estimated models are just too large. In comparison, svsCp does well in both cases. Overall it does the best in terms of selecting variables by maintaining low FDR and TotalMiss values. It also maintains good performance values.

TABLE 1
*Breiman simulation: $m = 400$, $n = 800$ and $105$ nonzero $\beta_j$*

|       | $\rho = 0$ (uncorrelated $X$) | | | | | $\rho = 0.9$ (correlated $X$) | | | | |
|-------|---|---|---|---|---|---|---|---|---|---|
|       | $\widehat{m}$ | $\text{pe}(\widehat{\mu})$ | TotalMiss | FDR | FNR | $\widehat{m}$ | $\text{pe}(\widehat{\mu})$ | TotalMiss | FDR | FNR |
| LARS  | 210.69 | 0.907 | 126.63 | 0.547 | 0.055 | 99.51  | 0.962 | 75.77  | 0.347 | 0.135 |
| svsCp | 126.66 | 0.887 | 61.14  | 0.323 | 0.072 | 58.86  | 0.952 | 66.38  | 0.153 | 0.164 |
| svsBMA| 400.00 | 0.918 | 295.00 | 0.737 | 0.000 | 400.00 | 0.966 | 295.00 | 0.737 | 0.000 |
| Step  | 135.53 | 0.876 | 70.35  | 0.367 | 0.075 | 129.24 | 0.884 | 137.10 | 0.552 | 0.208 |

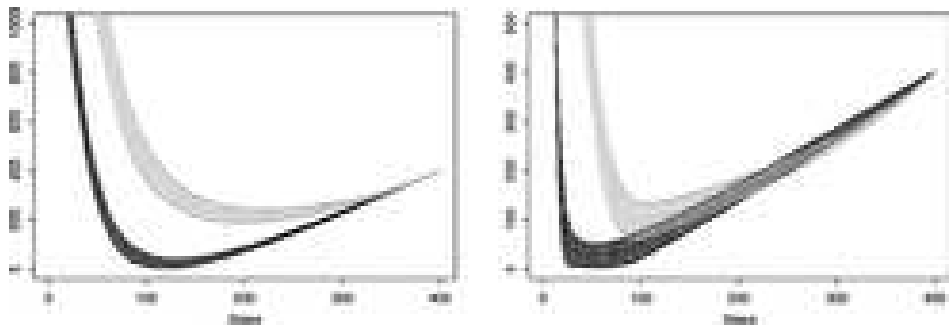

FIG. 3. *$C_p$ values from simulations where $\rho = 0$ (left) and $\rho = 0.9$ (right): bottom curves are from SVS; top curves are from LARS. The lines seen on each curve are the mean $C_p$ values based on 100 simulations. Note how the minimum value for SVS is near zero in both cases. Also superimposed on each curve are error bars representing mean values plus or minus one standard deviation.*

2. LARS's performance values are good, second only to svsBMA. However, low prediction error does not necessarily imply good variable selection.

**3. LARS $C_p$ values in orthogonal models.** Figure 3 shows the $C_p$ values for LARS from the two sets of simulations. It is immediately apparent that the $C_p$ curve in the uncorrelated case is too flat, leading to models which are too large. These simulations were designed to reflect an orthogonal design setting (at least asymptotically), so what is it about the orthogonal case that is adversely affecting LARS?

We can use Lemma 1 of the paper to gain some insight into this. For this argument I will assume that $m$ is fixed (the lemma is stated for $m = n$ but applies in general) and I will need to assume that $X_{n \times m}$ is a random orthogonal matrix, chosen so that its rows are exchangeable. To produce such an $X$, choose $m$ values $\mathbf{e}_{i_1}, \ldots, \mathbf{e}_{i_m}$ without replacement from $\{\mathbf{e}_1, \ldots, \mathbf{e}_n\}$, where $\mathbf{e}_j$ is defined as in Section 4.1, and set $X = [\mathbf{e}_{i_1}, \ldots, \mathbf{e}_{i_m}]$. It is easy to see that this ensures row-exchangeability. Hence, $\mu_1, \ldots, \mu_n$ are exchangeable and, therefore, $Y_i = \mu_i + \varepsilon_i$ are exchangeable since $\varepsilon_i$ are i.i.d. I will assume, as in (4.1), that $\varepsilon_i$ are independent $N(0, \sigma^2)$ variables.

For simplicity take $\sigma^2 = \overline{\sigma}^2 = 1$. Let $V_j$, for $j = 0, \ldots, m-1$, denote the $(j+1)$st largest value from the set of values $\{|Y_{i_1}|, \ldots, |Y_{i_m}|\}$. Let $k_0$ denote the true dimension, that is, the number of nonzero coordinates of the true $\boldsymbol{\beta}$, and suppose that $k$ is some dimension larger than $k_0$ such that $1 \leq k_0 < k \leq m \leq n$. Notice that $V_k \leq V_{k_0}$, and thus, by Lemma 1 and (4.10),

$$C_p(\widehat{\boldsymbol{\mu}}_k) - C_p(\widehat{\boldsymbol{\mu}}_{k_0}) = (V_k^2 - V_{k_0}^2) \sum_{j=1}^m \mathbb{1}\{|Y_{i_j}| > V_{k_0}\} + V_k^2 \sum_{j=1}^m \mathbb{1}\{V_k < |Y_{i_j}| \leq V_{k_0}\}$$

$$- \sum_{j=1}^m Y_{i_j}^2 \mathbb{1}\{V_k < |Y_{i_j}| \leq V_{k_0}\} + 2(k - k_0)$$

$$\leq -\Delta_k B_k + 2(k - k_0),$$

where $\Delta_k = V_{k_0}^2 - V_k^2 \geq 0$ and $B_k = \sum_{j=1}^m \mathbb{1}\{|Y_{i_j}| > V_{k_0}\}$. Observe that by exchangeability $B_k$ is a Binomial$(m, k_0/m)$ random variable. It is a little messy to work out the distribution for $\Delta_k$ explicitly. However, it is not hard to see that $\Delta_k$ can be reasonably large with high probability. Now if $k_0 > k - k_0$ and $k_0$ is large, then $B_k$, which has a mean of $k_0$, will become the dominant term in $\Delta_k B_k$ and $\Delta_k B_k$ will become larger than $2(k - k_0)$ with high probability. This suggests, at least in this setting, that $C_p$ will overfit if the dimension of the problem is high. In this case there will be too much improvement in the residual sums of squares when moving from $k_0$ to $k$ because of the nonvanishing difference between the squared order statistics $V_{k_0}^2$ and $V_k^2$.

**4. Summary.** The use of $C_p$ seems to encourage large models in LARS, especially in high-dimensional orthogonal problems, and can adversely affect variable selection performance. It can also be unreliable when used with stepwise regression. The use of $C_p$ with SVS, however, seems better motivated due to the benefits of model averaging, which mitigates the selection bias effect. This suggests that $C_p$ can be used effectively if model uncertainty is accounted for. This might be one remedy. Another remedy would be simply to use a different model selection criteria when using LARS.

Department of Biostatistics/Wb4
Cleveland Clinic Foundation
9500 Euclid Avenue
Cleveland, Ohio 44195
USA
e-mail: ishwaran@bio.ri.ccf.org